\documentclass{amsart}

\usepackage{hyperref}
\usepackage[utf8]{inputenc}
\usepackage[english]{babel}
\usepackage{amssymb}
\usepackage[shortalphabetic]{amsrefs}
\usepackage{etoolbox}
\usepackage{braket}
\usepackage{faktor}
\usepackage{mathtools}
\usepackage{stackrel}
\usepackage{verbatim}
\usepackage[english]{cleveref}

\newtheorem{theorem}{Theorem}[section]
\newtheorem{lemma}[theorem]{Lemma}
\newtheorem{corollary}[theorem]{Corollary}

\theoremstyle{definition}
\newtheorem{definition}[theorem]{Definition}

\theoremstyle{remark}
\newtheorem{remark}[theorem]{Remark}
\newtheorem{example}[theorem]{Example}

\newcommand*\Z{\mathbb{Z}}
\newcommand*\N{\mathbb{N}}
\newcommand*\Q{\mathbb{Q}}
\newcommand*\C{\mathbb{C}}

\newcommand*\Sn{\mathfrak{S}_n}
\newcommand{\defeq}{\stackrel{\textnormal{def}}{=}}
\newcommand{\Torg}{\operatorname{Tor}(\Sigma_g)}

\newcommand{\tdd}{\tilde{\dd}}
\newcommand{\tH}{\tilde{H}}
\newcommand{\Lo}{L_\omega}
\newcommand{\SSS}{\operatorname{S}}
\newcommand{\SG}{\mathfrak{S}}
\newcommand{\LLL}{\operatorname{\Lambda}}
\newcommand{\oomega}{\overline{\omega}}
\newcommand{\slg}{\mathfrak{sl}(2g)}
\newcommand{\spg}{\mathfrak{sp}(2g)}
\newcommand{\SV}[2]{\LLL^{#1} V \otimes \SSS^{#2} V}
\newcommand{\V}[2]{V_{#1 \omega_1 + \omega_{#2}}}
\newcommand{\W}[2]{W_{#1 \omega_1 + \omega_{#2}}}
\newcommand{\Uconf}[1]{\mathcal{C}_{#1}}
\newcommand{\Conf}[1]{\mathcal{F}_{#1}}
\newcommand*\Uconfn{\Uconf{n}}
\newcommand*\Confn{\Conf{n}}

\renewcommand{\Im}{\operatorname{im}}

\DeclareMathOperator{\dd}{d}
\DeclareMathOperator{\coker}{coker}
\DeclareMathOperator{\gr}{gr}

\newcommand{\ie}{i.e.\ }

\makeatletter
\newcommand*{\bigcdot}{%
  {\mathbin{\mathpalette\bigcdot@{}}}%
}
\newcommand*{\bigcdot@scalefactor}{.75}
\newcommand*{\bigcdot@widthfactor}{1.4}
\newcommand*{\bigcdot@}[2]{%
  \sbox0{$#1\vcenter{}$}
  \sbox2{$#1\cdot\m@th$}%
  \hbox to \bigcdot@widthfactor\wd2{%
    \hfil
    \raise\ht0\hbox{%
      \scalebox{\bigcdot@scalefactor}{%
        \lower\ht0\hbox{$#1\bullet\m@th$}%
      }%
    }%
    \hfil
  }%
}
\makeatother

\robustify{\bigcdot}

\begin{document}

\title[Cohomology of configuration spaces of surfaces]{Extra structure on the cohomology of configuration spaces of closed orientable surfaces}
\author{Roberto Pagaria}
\email{roberto.pagaria@unibo.it} 
\address{Dipertimento di matematica \\ Università di Bologna\\ Piazza di Porta San Donato 5\\ 40126 Bologna
\\Italy }

\begin{abstract}
The rational homology of unordered configuration spaces of points on any surface was studied by Drummond-Cole and Knudsen.
We compute the rational cohomology of configuration spaces on a closed orientable surface, keeping track of the mixed Hodge numbers and the action of 
the symplectic group on the cohomology.
We find a series with coefficients in the Grothendieck ring of $\spg$ that describes explicitly the decomposition of the cohomology into irreducible representations.
From that we deduce the mixed Hodge numbers and the Betti numbers, obtaining a new formula without cancellations.
\end{abstract}

\maketitle

\section{Introduction}
The ordered configuration space of $n$ points in a complex algebraic variety $X$ is
 \[\Confn(X)=\{(p_1, \dots ,p_n) \in X^n \mid p_i \neq p_j \}.\] 
We are interested in the unordered configuration space of $X$, that is
 \[\Uconfn(X)=\{I \subset X \mid |I|=n \} = \faktor{\Confn(X)}{\Sn}.\] 
We compute the rational cohomology of $\Uconfn(\Sigma_g)$ where $\Sigma_g$ is a Riemann surface of genus $g$.
Our computation is dual to the one by Drummond-Cole and Knudsen \cite{DCK}: they used the Chevalley–Eilenberg complex to compute $H_\bullet (\Uconfn(S);\Q)$ for any topological surface $S$ of finite type.
Since the mixed Hodge structure is defined for any algebraic varieties, there is no evident motivation for which the manipulations in \cite{DCK} and in \cite{Knudsen} are compatible with the Hodge structure.
Our work is based on the previous one by Félix and Tanré \cite{FT05} that uses the Cohen-Taylor spectral sequence to study the cohomology.

The Cohen-Taylor spectral sequence is a spectral sequence $E_\bigcdot (X,n)$ that converges to the rational cohomology of $\Confn(X)$, as proven in \cite{CohenTaylor78}*{pp.~117, 118}.
Kri\v{z} \cite{Kriz94} and Totaro \cite{Totaro96} used the Fulton and MacPherson's compactification \cite{FMacP94} to prove that for any $X$ smooth projective variety the spectral sequence degenerates at the second page, \ie $H(E_1^{\bigcdot,\bigcdot}(X,n), \dd_1)= \gr^W_\bigcdot H^\bigcdot (\Confn(X))$, where $W$ is the weight filtration defined by Deligne \cite{DelUtile}.

The symmetric group $\SG_n$ acts on $E_1(X,n)$ and the $\Sn$-invariant subalgebra computes the cohomology of the space $\Uconf{n}(X)$, indeed 
 \[H(E_1^{\bigcdot,\bigcdot}(X,n)^{\SG_n}, \dd_1) \simeq \gr^W_{\bigcdot} H^\bigcdot (\Uconf{n}(X)).\]
Furthermore, this isomorphism holds for all closed oriented manifolds, see \cite{FT05}*{Theorem 2}.
There exists an isomorphism $\gr^W_{\bigcdot} H^\bigcdot (\Uconf{n}(X)) \simeq H^\bigcdot (\Uconf{n}(X))$ both as algebras and as mixed Hodge structures, but in the case of Riemann surfaces $X=\Sigma_g$, this is not compatible with the action of the mapping class group.

Félix and Tanré \cite{FT05} presented the differential graded algebra $E_1^{\bigcdot,\bigcdot}(X,n)^{\SG_n}$ as a bigraded vector space $C_n$ with a differential $\dd$ and a complex multiplication law $\bigcdot$ that depends only on the cup product $\smile$ of $H^\bigcdot (X)$.
See \cite{FT05} for the definition in the general case or see \Cref{def:Cn} for the case $X=\Sigma_g$.

Using the result of Félix and Tanré, we construct an algebra $B_g$ with a filtration $\{F_n B_g\}$ and surjections of differential algebras $\varphi_n \colon B_g \twoheadrightarrow C_n$ that restricts to isomorphisms $F_nB_g \xrightarrow{\sim} C_n$.
The advantage of this method is that the algebra structure on $B_g$ is easy since it is an exterior algebra and the filtration $F_nB$ is induced by another grading of $B_g$.
On the other hand, $F_nB_g$ is not an algebra, so information about the ring structure of $ H^\bigcdot (\Uconf{n}(\Sigma_g))$ is lost.

We write $B$ as a shorthand for the algebra $B_g$.
We find an acyclic ideal $I$ of $B$, and then we define $A$ as the quotient $B/I$.
The algebra $A$ is filtered by $F_n A := \Im F_n B$, the induced filtration from $B$.
This filtration $F_nA$ is strictly compatible with the differential,  and this allow us to simplify the differential $\dd$ of $A$.
Finally, we obtain $H^\bigcdot (F_nA) \simeq \gr^W_\bigcdot H^\bigcdot ( \Uconfn( \Sigma_g))$ as algebras and as representations of the symplectic group.

Although the action of the mapping class group $\Gamma_g=\operatorname{MCG}(\Sigma_g)$ on $H^\bullet(\Uconf{n}(\Sigma_g))$ is not symplectic, it preserves the weight filtration $W$.
Hence the induced action on $\gr_\bigcdot^W H^\bigcdot ( \Uconf{n}(\Sigma_g))$ is symplectic and we study it as a representation of the symplectic group (see also Remark \ref{remark:symplectic}).

The next step is the explicit computation of the cohomology using the action of the Lie algebra $\spg$ on the model $(A,\dd)$.
From this analysis, we find out a formal power series with coefficients in $R_g$, \ie the Grothendieck ring of $\spg$.
For $g>0$ the following equation in $R_g[[t,s,u]]$ is proved in \Cref{theorem:main_theorem}:
\begin{multline}\label{eq:main_in_introduction}
\sum_{i,j,n} [\gr_{i+2j}^W H^{i+j}(\Uconf{n}(\Sigma_g))] t^i s^j u^n = \\
= \frac{1}{1-u}\Bigl( (1+t^2su^3)(1 + t^2u) + (1+t^2su^2) t^{2g}su^{2(g+1)} + (1+t^2su^2) \cdot \\
\cdot (1+t^2su^3)\sum_{\substack{1\leq j \leq g\\ i\geq 0}} [\V{i}{j}] t^{j+i} s^{i} u^{j+2i} (1+t^{2(g-j)}su^{2(g-j+1)}) \Bigr).
\end{multline}

Eq.~\eqref{eq:main_in_introduction} describes explicitly the decomposition of the associated graded module $\gr^W H^{\bullet}(\Uconf{n}(\Sigma_g))$ into irreducible representations.
Moreover, by taking the dimension $\dim \colon R_g \to \Z$, we obtain the mixed Poincaré polynomial of $\Uconfn(\Sigma_g)$ as the coefficient of $u^n$ in eq.~\eqref{eq:main_in_introduction}.
The dimension of the representations involved in our formula is calculated in \Cref{lemma:dim_rep}:
 \[\dim \V{i}{j}= \binom{2g+i+1}{i, j} \frac{2g+2-2j}{2g+2+i-j} \frac{j}{i+j}.\] 

The formula for the Betti numbers given in \cite{DCK} is different from the one in this paper, which has no cancellations and a more geometric meaning, because each summand corresponds to a specific submodule of the cohomology with a weight and a description of the symplectic group action.

We can give some explicit information about the mixed Hodge numbers.
\begin{corollary}
For $g>1$ the weights $h$ that appear in $H^k (\Uconfn(\Sigma_g))$ are in the range $\frac{3}{2}k-g-1 \leq h \leq \frac{3}{2}k$.
Moreover,  in this range the dimension of $\gr^W_h H^k(\Uconfn(\Sigma_g))$ is polynomial in $k$ of degree $2g-1$.
\end{corollary}
The polynomial growth of the Betti numbers was already established in \cite{DCK}*{Corollary 4.9}.

The cases of genera 0,1 have already been studied in \cites{Sevryuk,Salvatore04,Schiessl18} and in
\cites{Schiessl16,Maguire,Pagaria19}, respectively.
The Euler characteristic of the configuration spaces of any even-dimensional orientable closed manifold $M$ was computed by Félix and Thomas in \cite{FT00} and it is given by  the formula:
 \[ \sum_{n=0}^\infty \chi (\Uconfn(M)) u^n = (1+u)^{\chi (M)}.\] 
In the case of surfaces, this formula can be obtained from eq.~\eqref{eq:main_in_introduction} by setting $t=s=-1$ and taking the dimension of the representations.

\subsection*{Acknowledgements}\label{ackref}
I would like to thank Andrea Maffei and Sabino Di Trani for the useful discussions about representation theory.

\section{Models for \texorpdfstring{$H^\bigcdot (\Uconf{n}) $}{H(Cn)}}

Given a graded vector space $V$ we denote by $sV$ the same vector space with the degree shifted by one, \ie for any $v\in V$ of degree $i$, the element $sv \in sV$ has degree $i+1$.
We fix a symplectic base of $H^\bigcdot (\Sigma_g)$ whose elements are $1 \in H^0(\Sigma_g)$, $a_1, \dots, a_g, b_1,\dots, b_g \in H^1(\Sigma_g)$, $p \in H^2(\Sigma_g)$.
The cup product $\smile$ in $H^\bigcdot (\Sigma_g)$ is given by $a_i \smile b_j = 0$ for $i \neq j$ and $a_i \smile b_i =p$ for all $i=1, \dots, g$.

\begin{definition}
\label{def:Cn}
Let $C_n$ be the bigraded vector space
 \[ C_n \defeq  \bigoplus_{r=0}^{\lfloor \frac{n}{2} \rfloor} \LLL^{n-2r} \left( H(\Sigma_g) \right) \otimes \LLL^{r} \left( sH(\Sigma_g) \right) \] 
where the bigrade of $H^k(\Sigma_g)$ is $(k,0)$, the bigrade of $sH^k(\Sigma_g)$ is $(k,1)$, and the graded-symmetric algebra is constructed with respect to the total degree.
We endow $C_n$ with the product:
\begin{multline*}
(x_1 \wedge \dots \wedge x_{n-2r} \otimes sy_1 \wedge \dots sy_r) \bigcdot (z_1 \wedge \dots \wedge z_{n-2t} \otimes sw_1 \wedge \dots sw_t) \defeq \\ \frac{1}{(n-2r-2t)!}\sum_{\substack{\sigma \in \SG_{n-2r} \\ \tau \in \SG_{n-2t}}} \pm \alpha_{\sigma,\tau} \otimes (\beta_{\sigma} \wedge \gamma_{\tau}),
\end{multline*}
where:
\begin{align*}
& \alpha_{\sigma,\tau}= (x_{\sigma(1)} \smile z_{\tau(1)}) \wedge \dots \wedge (x_{\sigma(n-2r-2t)} \smile z_{\tau(n-2r-2t)}) \\
& \beta_{\sigma}= s(x_{\sigma(n-2r-2t+1)} \smile x_{\sigma(n-2r-2t+2)}\smile w_1) \wedge \dots \wedge s(x_{\sigma(n-2r-1)} \smile x_{\sigma(n-2r)}\smile w_t) \\
& \gamma_{\tau}= s(z_{\tau(n-2r-2t+1)} \smile z_{\tau(n-2r-2t+2)}\smile y_1) \wedge \dots \wedge s(z_{\tau(n-2t-1)} \smile z_{\tau(n-2t)}\smile y_r),
\end{align*}
and the sign $\pm$ is given by the Koszul rule.
We consider the differential 
 \[\dd \colon \LLL^{n-2r} \left( H(\Sigma_g) \right) \otimes \LLL^{r} \left( sH(\Sigma_g) \right) \to \LLL^{n-2r+2} \left( H(\Sigma_g) \right) \otimes \LLL^{r-1} \left( sH(\Sigma_g) \right)\] 
of degree $(2,-1)$ defined on the generators by 
\begin{multline*}
\dd(1 \wedge \dots \wedge 1 \otimes sy) \defeq \frac{1}{2} \Bigl(1 \wedge \dots \wedge 1 \wedge y \wedge p + 1 \wedge \dots \wedge 1  \wedge (y \smile p) \\ - \sum_{i=1}^g 1 \wedge \dots \wedge 1 \wedge (y \smile a_i) \wedge b_i + \sum_{i=1}^g 1 \wedge \dots \wedge 1 \wedge (y \smile b_i) \wedge a_i \Bigr).
\end{multline*}
\end{definition}

For the sake of notation, in the following we will write $x_1 \wedge \dots \wedge x_{k} \otimes sy_1 \wedge \dots sy_r$ for the element $x_1 \wedge \dots \wedge x_{k} \wedge 1\wedge \dots \wedge 1 \otimes sy_1 \wedge \dots sy_r$ where the number of $1$ omitted is $n-2r-k$.

\begin{theorem}[\cite{FT05}*{Theorems 1, 14}]
The triple $(C_n, \bigcdot, \dd)$ is a differential graded algebra and it is isomorphic to the $\SG_n$-invariants of the first page $E_1$ of the Cohen-Taylor spectral sequence for $\Uconfn(\Sigma_g)$:
 \[ (E_1,\dd_1)^{\SG_n} \cong  (C_n, \dd).\] 
\end{theorem}

\begin{definition}
\label{def:Bg}
Let $B=B_g= \LLL^\bigcdot (\tH^\bigcdot(\Sigma_g)\oplus sH^\bigcdot (\Sigma_g))$ be the graded-symmetric algebra on the trigraded vector space $\tH^\bigcdot(\Sigma_g)\oplus sH^\bigcdot (\Sigma_g)$, where the grading is given in \Cref{tab:deg} and the total degree is $|x|=\deg_1(x)+\deg_2(x)$.
\begin{table}
\centering
\begin{tabular}{c|c|c|c}
 & $\deg_1$ & $\deg_2$ & $\deg_3$ \\ 
\hline 
$a_i, b_i$ & 1 & 0 & 1 \\ 
\hline 
$p$ & 2 & 0 & 1 \\ 
\hline 
$s1$ & 0 & 1 & 2 \\ 
\hline 
$sa_i, sb_i$ & 1 & 1 & 2 \\ 
\hline 
$sp$ & 2 & 1 & 2 \\ 
\end{tabular} 
\caption{The degree of the generators of $B_g$.}\label{tab:deg}
\end{table}
We endow $B_g$ with the following differential:
\begin{equation}
\label{eq:def_diff_B}
\begin{aligned}
&\dd(x)=0 \quad \textnormal{for } x \in \tH^\bigcdot(\Sigma_g) \\
&\dd(sp)= p^2, \quad \dd(s1)= p- \sum_{i=1}^g a_i b_i \\
&\dd(sa_i)= a_i p, \quad \dd(sb_i)= b_i p \quad \textnormal{for } i=1, \dots, g.
\end{aligned}
\end{equation}
\end{definition}

\begin{remark}
The elements of $C_n$ (resp.\ of $B$) can be interpreted geometrically as follows.
Generators $a_i$ and $b_i$ for $i=1, \dots, g$ are an average over all particles (\ie points of the configuration) of the motion of that particle along the curve $a_i$ (resp.\ $b_i$).
The element $p$ is the average over all particles of the motion of that particle on the entire surface.
The generator $s1$ is the average over all pairs of particles of the rotation of one particle around the other. 
Similar description holds for the other generators.

In order to describe cohomological classes we need to resolve the collision problems of the moving particle with the other particles in the configuration.  This is possible only if its differential is zero. 

The only difference between the generators in $B$ and in $C_n$ consists in the multiplication by a numerical coefficients, as shown in the following Lemma \ref{lemma:hom_phi_n}.
\end{remark}

Notice that the differential $\dd$ is compatible with $\deg_1$ and $\deg_2$ of degree $(2,-1)$, but not with the third grading since $\dd(s1)=p - \sum_{i=1}^g a_ib_i$.
For all $n\in \N$ consider the morphisms $\varphi_n \colon B \to C_n$ defined by
\begin{align*}
& \varphi_n(a_i) = a_i, \quad  \varphi_n(sa_i) = (n-1)!^2(n-1)(n-1+g)sa_i \quad \textnormal{for } i=1, \dots g,\\
& \varphi_n(b_i) = b_i, \quad  \varphi_n(sb_i) = (n-1)!^2(n-1)(n-1+g) sb_i \quad \textnormal{for } i=1, \dots g, \\
&\varphi_n (p)= (n-1)!(n-1+g)p, \quad \varphi_n (sp)= 2(n-1)!^3(n-1)(n-1+g)^2sp, \\
&\varphi_n(s1) = (n-1)!(n-1) s1.
\end{align*}

\begin{lemma}\label{lemma:hom_phi_n}
The map $\varphi_n \colon B \to C_n$ is a morphism of differential graded algebras.
\end{lemma}
\begin{proof}
We verify the compatibility between $\varphi_n$ and the differential.
We first compute products in $C_n$:
 \[
a_i \bigcdot b_i = \frac{1}{n!} \sum_{\substack{\sigma \in \SG_n \\ \tau \in \SG_n}} \alpha_{\sigma,\tau} = (n-1)!(n-1) a_i \wedge b_i + (n-1)! p,\] 
analogously we have
\begin{gather*}
a_i \bigcdot p = (n-1)!(n-1) a_i \wedge p, \\
b_i \bigcdot p = (n-1)!(n-1) b_i \wedge p, \\
p \bigcdot p = (n-1)!(n-1) p \wedge p.
\end{gather*}

Now we prove the equalities $\dd (\varphi_n(x))= \varphi_n(\dd(x))$ for all the generators $x$:
for $x=s1$ we have
\begin{align*}
\dd (\varphi_n(s1)) &= \dd ((n-1)!(n-1) s1) \\
&= (n-1)!(n-1) (p- \sum_{i=1}^g a_i \wedge b_i) \\
&= (n-1)!(n-1) p - \sum_{i=1}^g a_i \bigcdot b_i - (n-1)! p \\
&= (n-1)!(n-1+g) p  - \sum_{i=1}^g a_i \bigcdot b_i \\
&= \varphi_n(p- \sum_{i=1}^g a_i b_i) = \varphi_n(\dd(s1)),\\
\end{align*}
for $x=sp$
\begin{align*}
\dd (\varphi_n(sp)) &= \dd (2(n-1)!^3(n-1)(n-1+g)^2 sp) \\
&=2(n-1)!^3(n-1)(n-1+g)^2 \frac{1}{2}p\wedge p \\
&= (n-1)!^2(n-1+g)^2 p \bigcdot p \\
&= \varphi_n(p^2) = \varphi_n(\dd(sp)),\\
\end{align*}
for $x=sa_i$, $i=1,\dots,g$
\begin{align*}
\dd(\varphi_n(sa_i)) &= \dd ((n-1)!^2(n-1)(n-1+g)sa_i) \\
&= (n-1)!^2(n-1)(n-1+g) a_i \wedge p\\
&= (n-1)!(n-1+g)  a_i \bigcdot p \\
&= \varphi_n(a_i p) = \varphi_n(\dd(sa_i)),
\end{align*}

and analogously for $sb_i$.
For all other generators $x= a_i, b_i,p$ we have $\dd(\varphi_n(x))=0=\varphi_n(\dd(x))$.
\end{proof}

\begin{definition}
Let $\left\{ F_k B \right\}_{k\in \N}$ be the filtration of $B$ defined by the third grading, i.e\ $F_kB = \oplus_{i\leq k} B^{\bigcdot, \bigcdot, i}$.
\end{definition}
Notice that the inclusion $\dd (F_kB) \subseteq F_kB$ holds since it holds for all generators of $B$.

\begin{lemma}\label{lemma:sur_phi_n}
Let $n-1+g \neq 0$ and $k+2r \leq n$.
For all $x_1, \dots x_k, y_1, \dots y_r \in H(\Sigma_g)$ we have in $C_n$ the following equality
\begin{multline*}  1 \wedge \dots \wedge 1 \wedge x_1 \wedge \dots \wedge x_k \otimes sy_1 \wedge \dots \wedge sy_r = \\
= \lambda \varphi_n(x_1) \bigcdot \dots \bigcdot \varphi_n(x_k) \bigcdot \varphi_n(sy_1) \bigcdot \dots \bigcdot \varphi_n(sy_r) +z
\end{multline*}
for some $\lambda \in \Q^*$ and some $z \in \varphi_n(F_{k+2r-1}B)$.
In particular $\varphi_{n} \colon F_nB \to C_n$ is surjective.
\end{lemma}
\begin{proof}
We first consider the case $k=0$: we prove the statement by induction on $r$.
The base step for $r=1$ follows from the definition of $\varphi_n$ and from $n-1+g \neq 0$.
Suppose $r\geq 2$ and consider the product $sy_1 \wedge \dots \wedge sy_{r-1} \bigcdot sy_r$ in $C_n$:
\begin{align*}
sy_1 \wedge \dots \wedge sy_{r-1} \bigcdot sy_r &= \frac{1}{(n-2r)!} \sum_{\substack{\sigma \in \SG_{n-2r+2}\\ \tau \in \SG_{n-2}}} \beta_{\sigma,\tau} \wedge \gamma_{\sigma,\tau} \\
&= \frac{(n-2r+2)!(n-2)!}{(n-2r)!}sy_1 \wedge \dots \wedge sy_{r-1} \wedge sy_r.
\end{align*}

By inductive hypothesis we have $sy_1 \wedge \dots \wedge sy_{r-1}= \lambda' \varphi_n(sy_1) \bigcdot \dots \bigcdot \varphi_n(sy_{r-1})$, so $sy_1 \wedge \dots \wedge sy_r 
= \lambda \varphi_n(sy_1) \bigcdot \dots \bigcdot \varphi_n(sy_r)$ for some $\lambda \in \Q^*$.

Now we proceed by induction on $k$, the base step $k=0$ is already been proved.
Suppose $k>0$, we have
\begin{gather*}
1 \wedge \dots \wedge 1 \wedge x_1 \wedge \dots \wedge x_{k-1} \otimes sy_1 \wedge \dots \wedge sy_r \bigcdot x_k = \frac{1}{(n-2r)!	} \sum_{\substack{\sigma \in \SG_{n-2r} \\ \tau \in \SG_{n}}}\alpha_{\sigma,\tau} \otimes \beta_{\sigma,\tau} \\
\begin{split}
	=\, & (n-1)! \sum_{i=1}^{k-1} x_1 \wedge \dots \wedge (x_i \cup x_k) \wedge \dots \wedge x_{k-1} \otimes sy_1 \wedge \dots \wedge sy_r \\
	&+ 2(n-1)! \sum_{j=1}^{r} x_1 \wedge \dots  \wedge x_{k-1} \otimes sy_1 \wedge \dots \wedge s(y_j \cup x_k) \wedge \dots \wedge sy_r \\
	&+ (n-1)! (n-k-2r+1)  x_1 \wedge \dots \wedge x_k \otimes sy_1 \wedge \dots \wedge sy_r.
\end{split}
\end{gather*}

The first two sums belong to $\varphi_n(F_{k+2r-1}B)$ by inductive hypothesis.
We also have
 \[ x_1 \wedge \dots \wedge x_{k-1} \otimes sy_1 \wedge \dots \wedge sy_r = \lambda'' \varphi_n(x_1) \bigcdot \dots \bigcdot \varphi_n(x_{k-1}) \bigcdot \varphi_n(sy_1) \bigcdot \dots \bigcdot \varphi_n(sy_r) + z\] 
for some $\lambda'' \in \Q^*$ and $z \in \varphi_n(F_{k+2r-2}B)$.
Using that $x_k$ is a multiple of $\varphi_n(x_k)$ and that $z \bigcdot \varphi_n(x_k) \in \varphi_n(F_{k+2r-1}B)$, we obtain the claimed equality.
\end{proof}

\begin{lemma}\label{lemma:iso_phi_n}
The restricted chain map $\varphi_n \colon F_nB \to C_n$ is an isomorphism.
\end{lemma}
\begin{proof}
\Cref{lemma:hom_phi_n} ensures that $\varphi_n$ is a homomorphism of chain complexes and \Cref{lemma:sur_phi_n} gives the surjectivity.
We complete the proof with a dimensional argument.

Let $S(t)$ be the formal power series $\frac{(1+t)^{2g}}{(1-t)^2}$ and $R(t)=\frac{(1+t)^2}{(1-t)^{2g}}$, they are the Poincaré series of $\LLL^\bigcdot (H(\Sigma_g))$ and of $\LLL^\bigcdot (sH(\Sigma_g))$ with respect to the cohomological degree.
For any power series $K(t)$, we denote the coefficient of $t^k$ by $[t^k]K(t)$.
The dimension of $C_n$ is
 \[ \dim C_n = \sum_{r=0}^{\lfloor \frac{n}{2} \rfloor} [t^{n-2r}]S(t) \cdot [t^r]R(t) = [t^n](S(t)R(t^2)).\] 
By definition of $F_nB$, we have
 \[ \dim F_nB = \sum_{i=0}^n \dim B^{\bigcdot, \bigcdot, i}= \sum_{i=0}^n [t^i] \left( \frac{(1+t)^{2g}(1+t^2)^{2}}{(1-t)(1-t^2)^{2g}} \right). \] 
Let $K(t)$be the power series $ \frac{(1+t)^{2g}(1+t^2)^{2}}{(1-t)(1-t^2)^{2g}} $.
The equality
\begin{align*}
[t^n](S(t)R(t^2)) &= [t^n] \left( \frac{(1+t)^{2g}(1+t^2)^{2}}{(1-t)^2(1-t^2)^{2g}} \right)\\
 &= \sum_{i=1}^n [t^{n-i}]\left(\frac{1}{1-t}\right) [t^i]K(t)\\
 &= \sum_{i=1}^n [t^i]K(t)
\end{align*}
completes the proof.
\end{proof}

Consider the ideal $I$ of $B$ generated by $sp$ and $p^2$ and define $A$ as the quotient $B/I$.
Since $\dd(I) \subseteq I$, $A$ is a differential graded algebra.
The filtration $\{F_nB \}$ induces two filtrations $\{F_nI \}$ and $\{F_nA \}$ on $I$ and $A$, respectively.
\begin{lemma}\label{lemma:I_is_acyclic}
The ideal $I$ is acyclic, \ie $H(I,\dd)=0$.
Moreover the chain complexes $(F_n I, \dd)$ are acyclic for all $n$, \ie $H(F_n I,\dd)=0$.
\end{lemma}
\begin{proof}
First, notice that the filtration $F_nI$ is induced by the third grading $\deg_3$ of $B$ and that $p$ and $sp$ are homogeneous elements of degrees $\deg_3(p)=1$ and $\deg_3(sp)=2$, hence:
\begin{equation} \label{eq:hom_I}
F_nI := F_nB \cap I = F_{n-2}p^2 + F_{n-2}sp.
\end{equation}

Since $B$ is an exterior algebra, we have that $(p^2)\cap (sp)= (p^2 sp)$ and, keeping track of the gradation,  $(p^2)\cap (sp) \cap F_nB = p^2 sp F_{n-4}B$.

Consider a generic element $xp^2+ y sp$ of $F_nI$ homogeneous with respect to $\deg_1$ and $\deg_2$, and suppose that $xp^2+ y sp$ belongs to $\ker \dd$, we will prove that $xp^2+ y sp \in \dd(F_nI)$.
Since $\dd$ preserves the first two degrees, we can assume $x$ and $y$ to be homogeneous with respect to $\deg_1$ and $\deg_2$.
By eq.~\eqref{eq:hom_I} we can suppose that $x$ and $y$ are in $F_{n-2}B$.
We use the hypothesis $\dd(xp^2+ y sp)=0$:
 \[ \dd(xp^2+ y sp)=\dd(x) p^2 + (-1)^{|y|}y p^2 + \dd(y) sp=0. \] 
The element $\dd(x) p^2 + (-1)^{|y|}y p^2$ belongs to $(p^2) \cap (sp) \cap F_nB$, so there exists $z \in F_{n-4}B$ such that $\dd(x) p^2 + (-1)^{|y|}y p^2 = z p^2 sp$.
We use again the fact that $B$ is an exterior algebra to obtain $\dd(x) + (-1)^{|y|}y = z sp$.
Therefore, we have
 \[ xp^2+ y sp =xp^2 - (-1)^{|y|}\dd(x) sp = p^2 x  - sp \dd(x) = \dd (sp\, x).\] 
Since $x \in F_{n-2}B$,  then $sp\, x \in F_n I$.
We have proven that $H(F_nI,\dd)=0$, the vanishing of $H(I,\dd)$ follows from $I = \cup_n F_nI$.
\end{proof}

Notice that $A$ is isomorphic to $\LLL^\bigcdot \left( \tH(\Sigma_g) \oplus sH^{\leq 1}(\Sigma_g) \right) /(p^2)$.

\section{Some facts of representation theory}
\label{sect:repr_theory}

Let $\Sigma_g$ be a Riemann surface, $\Gamma_g$ be its mapping class group, and $\Torg$ the Torelli subgroup.
Recall the short exact sequence
 \[ 0 \to \Torg \to \Gamma_g \to \operatorname{Sp}(2g; \Z) \to 0 .\] 

Consider the weight filtration $W_\bigcdot$ of the cohomology of the algebraic variety $\Uconf{n}(\Sigma_g)$, definition and properties of this filtration can be found in \cite{DelUtile}.
We define the module $\gr^W_i H^{\bigcdot}(\Uconf{n}(\Sigma_g))$ as the quotient $W_i H^{\bigcdot}(\Uconf{n}(\Sigma_g))/W_{i-1} H^{\bigcdot}(\Uconf{n}(\Sigma_g))$ and their direct sum $\gr_\bigcdot^W H^{\bigcdot}(\Uconf{n}(\Sigma_g)) =\oplus_{i} \gr^W_i H^{\bigcdot}(\Uconf{n}(\Sigma_g))$ is a bigraded ring.

The natural action of the subgroup $\Torg$ on $H^i(\Uconfn(\Sigma_g))$ may be non-trivial, but the induced action on $\gr^{W}_\bigcdot H^\bigcdot (\Uconfn(\Sigma_g))$ is trivial.
Indeed, since $\Sigma_g$ is a compact algebraic variety, Totaro proved in \cite{Totaro96}*{Theorem 3} that the Leray filtration for the inclusion $\Confn(\Sigma_g) \hookrightarrow \Sigma_g^n$ coincides with the weight filtration on $H^\bigcdot (\Confn(\Sigma_g))$.
Therefore, each homeomorphism of the pair $(\Confn(\Sigma_g),\Sigma_g^n)$ preserves the Leray filtration, and thus the weight filtration $W_\bigcdot$.
In particular this applies to each element in $\Gamma_g$ acting on the pair $(\Confn(\Sigma_g),\Sigma_g^n)$.
Indeed, $\gr^{W}_\bigcdot H^\bigcdot (\Uconfn(\Sigma_g))\simeq H (E_1^{\SG_n},\dd_1)$ functorially, hence the isomorphism is $\Gamma_g$-equivariant.
The action of $\Gamma_g$ on the algebra $E_1^{\SG_n}$ is clearly symplectic thus $\Torg$ acts trivially on $\gr^{W}_\bigcdot H^\bigcdot (\Uconfn(\Sigma_g))$.
The action of the Torelli group is studied in \cite{AB19} in the case of once punctured surfaces and it is non-trivial on $H_2(\Uconfn(\Sigma_g \setminus \{*\}))$; the case of compact surfaces is similar.

\begin{remark}\label{remark:symplectic}
From \Cref{theorem:main_theorem}, we deduced that the filtration $W_\bigcdot$ is trivial in cohomological degrees $0,1,2$ and also in degree $3$ if $g=2$, since the graded module $\gr_\bigcdot^W H^i(\Uconfn(\Sigma_g))$ is concentrated in a unique degree for $i=0,1,2$ (and $i=3,4$ if $g=2$).
Thus in these cases the action of the mapping class group is symplectic.

Looijenga in \cite{Looijenga} proves that the action of $\Torg$ on $H^3(\Uconf{3}(\Sigma_g))$ is non-trivial for $g\geq 3$ and so in these cases the action is not symplectic.
\end{remark}

We consider $\gr^{W}_\bigcdot H^\bigcdot (\Uconfn(\Sigma_g);\C)$ as a representation of the Lie algebra $\spg$ associated to the complex symplectic group.
If we denote the fundamental weights of $\spg$ by $\omega_1, \dots, \omega_g$, the irreducible representations of $\spg$ are the highest weight representation $V_\lambda$ for all dominant weights $\lambda=\sum_{i=1}^g \lambda_i \omega_i$, $\lambda_i \in \N$.
The cohomology of $\Sigma_g$ in degree one is given by the standard representation, \ie $H^1(\Sigma_g)=V_{\omega_1}$.

Let $V$ be the $\spg$-representation  $=H^1(\Sigma_g)\simeq V_{\omega_1}$ .
Before computing the cohomology of $(A,\dd)$ we need to know the cohomology of $(\SV{\bigcdot}{\bigcdot}, \tdd)$, where the differential $\tdd$ is defined by $\tdd (1 \otimes v)= v \otimes 1$ and $\tdd (v \otimes 1)=0$ for all $v \in V$.
The standard action of $\slg$ on $V$ induces an action on  $(\SV{\bigcdot}{\bigcdot}, \tdd)$, since the differential $\tdd$ is $\slg$-equivariant.

We will call $\omega_1, \dots, \omega_{2g-1}$ the fundamental weights of $\slg$ and $W_\mu$ its irreducible representations associated to a dominant weight $\mu=\sum_{i=1}^{2g-1} m_i \omega_i$, $m_i \in \N$.
\begin{lemma}\label{lemma:sl_dec}
The $\slg$-representation $\SV{j}{i}$ decomposes, for $j\leq 2g$, as 
 \[\SV{j}{i} =\W{i}{j} \oplus \W{(i-1)}{j+1}.\] 
\end{lemma}

\begin{proof}
It is known that $\SSS^i W_{\omega_1}=W_{i\omega_1}$ and $\LLL^j W_{\omega_1} = W_{\omega_j}$.
Let $\sigma=(1,j+1) \in \SG_{2g}$ be an element of the Weyl group of $\slg$.
The element $i\omega_1+\sigma(\omega_j)=(i-1) \omega_1+\omega_{j+1}$ is a dominant weight for $i>0$.
By the Parthasarathy–Ranga-Rao–Varadarajan conjecture (see \cites{Littelmann,Kumar88}) $\W{i}{j}$ and $\W{(i-1)}{j+1}$ are contained in the tensor product $W_{i\omega_1}\otimes W_{\omega_j}$.
Use the Weyl dimension formula to find
 \[ \dim \W{i}{j} = \binom{i+j-1}{i}\binom{i+2g}{i+j}.
\] 
The equality $\binom{i+j-1}{i}\binom{i+2g}{i+j}+\binom{i+j-1}{i-1}\binom{i+2g-1}{i+j}=\binom{2g}{j}\binom{i+2g-1}{i} $ completes the proof.
\end{proof}

\begin{lemma}
\label{lemma:exactness_tilde_dd}
The differential complex $(\SV{\bigcdot}{\bigcdot},\tilde{\dd})$ is exact in positive degree.
\end{lemma}

\begin{proof}
The differential 
 \[\tilde{\dd} \colon W_{i\omega_1 +\omega_j} \oplus W_{(i-1)\omega_1 +\omega_{j+1}} \to W_{(i-1)\omega_1 +\omega_{j+1}} \oplus W_{(i-2)\omega_1 +\omega_{j+2}}\] 
is a non-zero morphism of representations.
Therefore, we have $W_{i\omega_1 +\omega_j} = (\ker \tilde{\dd})^{j,i}$ and $W_{(i-1)\omega_1 +\omega_{j+1}}= (\Im \tilde{\dd})^{j+1,i-1}$ for $j\geq 0$.
Obviously $(\Im \tilde{\dd})^{0,0}=0$, so the equality
 \[ H^\bigcdot (\SV{\bigcdot}{\bigcdot},\tilde{\dd})= \C\] 
completes the proof.
\end{proof}

\begin{remark}
The complex $(\SV{\bigcdot}{\bigcdot}, \tilde{\dd})$ is the Koszul resolution of the trivial $\SSS^\bigcdot V$-module $\C$, hence it is an exact complex.
\end{remark}

Since the Lie algebra $\slg$ does not act on $A_g$, we need to present a branching rule for $\spg \subset \slg$.
For the sake of an uniform notation, we define $V_\lambda=0$ if $\lambda$ is not a dominant weight.
\begin{lemma}[Branching rule]
\label{lemma:branching}
The $\slg$-module $\W{i}{j}$ decomposes as $\spg$-module in the following ways:

\begin{align*}
&\W{i}{j} = \bigoplus_{k=0}^{\lfloor \frac{j-1}{2} \rfloor} \V{i}{j-2k} \oplus \bigoplus_{k=0}^{\lfloor \frac{j-2}{2} \rfloor} \V{(i-1)}{j-2k-1} \quad \textnormal{if } 2\leq j\leq g, \\
&\W{i}{j} = \bigoplus_{k=0}^{\lfloor \frac{2g-j-1}{2} \rfloor}\V{i}{2g-j-2k} \oplus \bigoplus_{k=0}^{\lfloor \frac{2g-j-2}{2} \rfloor} \V{(i-1)}{2g-j-2k-1} \quad \textnormal{if } j \geq g.
\end{align*}

\end{lemma}

\begin{proof}
We apply the result of \cite{SchumannTorres}*{Theorem 1}.
The diagram associated to $i\omega_1+\omega_j$ has a hook shape with row length $i+1$ and column length $j$.
Fill each box with labels in the ordered set $\{1<...<g<\bar{g}<...<\bar{1}\}$, such that it becomes a semi-standard Young tableau (SSYT) \ie the rows are non-decreasing and columns are increasing.
The word $w(T)$ -- associated to a SSYT $T$ -- is the word obtained by reading the tableaux from right to left and from top to bottom.
By convention, $e_{\bar{a}}=-e_a$.
A word $w(T)=a_1a_2 \dots a_k$ is admissible if for each $r\leq k$ the element $\sum_{s=1}^r e_{a_s}$ is a dominant weight for $\spg$.
The decomposition of $\W{i}{j}$ into $\spg$-representations is given by
 \[ \W{i}{j} = \bigoplus_{w(T) \textnormal{ admissible}} V_{\lambda(T)},\] 
where $\lambda(T)=\sum_{s=1}^{|w(T)|} e_{a_s}$.

Suppose $w(T)$ is admissible, then the first row of $T$ is labelled only by ones.
For $j\leq g$, all possible labels of the first column of $T$, from top to bottom, are the following:
\begin{itemize}
\item $1,2,\dots, j-k, \overline{j-k},\overline{j-k-1}, \dots, \overline{j-2k+1}$, where $k$ is an integer such that $0\leq 2k \leq j-1$ 
\item $1,2,\dots,j-k-1,\overline{j-k-1}, \dots, \overline{j-2k},\bar{1}$,
where $k$ is an integer such that $0\leq 2k \leq j-2$ and $i>0$.
\end{itemize}
Our decomposition follows, the case $j>g$ being analogous.
\end{proof}

Let $\omega$ be the element $\sum_{i=1}^{g} a_i \wedge b_i \in \Lambda^\bigcdot H^1(\Sigma_g) \subset A$.
The differential $\dd$ of $A$ involves the multiplication by $\omega$ (see eq.\ \eqref{eq:def_diff_B}).
Thus we need to study the operator $\Lo \colon \SV{j}{i} \to \SV{j+2}{i}$ defined by left multiplication by $\omega$.

\begin{lemma}[\cite{FHbook}*{Theorem 17.5}]
\label{lemma:dec_ext_power}
The $\spg$-representation $\LLL^j V$ is isomorphic to $\LLL^{2g-j} V$ and decomposes, for $j\leq g$, as
 \[\LLL^j V = W_{\omega_j} = \bigoplus_{k=0}^{\lfloor \frac{j}{2} \rfloor} V_{\omega_{j-2k}}.\] 
Moreover, $(\ker \Lo )^{2g-j}=V_{\omega_j} \subset \LLL^{2g-j} V$ and $(\coker \Lo)^j = V_{\omega_j} \subset \LLL^{j} V$.
\end{lemma}

We denote by $R_g$ the Grothendieck ring of $\spg$, \ie $R_g$ is the free $\Z$-module with basis the irreducible (finite dimensional) representations of $\spg$.
The ring structure on $R_g$ is induced by the tensor product of representations, however we do not need the multiplicative structure.

\begin{lemma}\label{lemma:dec_tensor_prod}
For $i>0$ and $1 \leq j\leq g$, we have
\begin{equation}\label{eq:dec_tensor_prod}
V_{\omega_j} \otimes \SSS^{i} V = \V{i}{j} \oplus \V{(i-1)}{j+1} \oplus \V{(i-1)}{j-1} \oplus \V{(i-2)}{j}.
\end{equation}

\end{lemma}

\begin{proof}
We use \Cref{lemma:dec_ext_power,lemma:sl_dec,lemma:branching}:
\begin{align*}
V_{\omega_j} \otimes \SSS^{i} V &= \SV{j}{i}\ominus \SV{j-2}{i} \\
&= W_{i\omega_1+\omega_j} \oplus W_{(i-1)\omega_1+\omega_{j+1}} \ominus W_{i\omega_1+\omega_{j-2}} \ominus W_{(i-1)\omega_1+\omega_{j-1}}\\
&=\V{i}{j} \oplus \V{(i-1)}{j+1} \oplus \V{(i-1)}{j-1} \oplus \V{(i-2)}{j},
\end{align*}
where the symbol $\ominus$ is the negation in the Grothendieck ring $R_g$.
\end{proof}

\begin{definition}
Let $W= \oplus_{i,j} W^{i,j}$ be a bigraded representation of the symplectic Lie algebra $\spg$. The Hilbert–Poincaré series of $W$ is the formal power series
 \[ P_W(t,s)= \sum_{i,j}[W^{i,j}] t^i s^j \quad \in R_g[[t,s]].\] 
\end{definition}

Recall that the bidegree of $v \otimes 1$ is $(1,0)$ and the one of $1 \otimes v$ is $(1,1)$, for all $v \in V$.

\begin{corollary}
The Hilbert–Poincaré series of the representation $\SV{\bigcdot}{\bigcdot}$ is
\begin{align*}
P_{\SV{}{}}(t,s)=& \frac{t^{2(g+1)}-1}{t^2-1}+t^2s \frac{t^{2g}-1}{t^2-1} +\\
&+ (1+s)(1+t^2s)\sum_{\substack{1 \leq j \leq g \\i\geq 0 }} [\V{i}{j}] \frac{t^{2(g-j+1)}-1}{t^2-1}t^{i+j}s^i.
\end{align*}
\end{corollary}
\begin{proof}
We first notice that \Cref{lemma:dec_ext_power} implies
 \[ \sum_{j=0}^{2g} [\LLL^j V] t^j = \sum_{j=0}^g \sum_{k=0}^{g-j} [V_{\omega_j}] t^{j+2k} = \sum_{j=0}^g [V_{\omega_j}] \frac{t^{2(g-j+1)}-1}{t^2-1}t^{j},\] 
where we set $[V_{\omega_0}]=1 \in R_g$.
For $1\leq j \leq g$, \Cref{lemma:dec_tensor_prod} implies that
\begin{multline*}
\sum_{i\geq 0} [V_{\omega_j} \otimes \SSS^i V] (ts)^{i} =  \sum_{i\geq 0} [V_{i \omega_1+\omega_j}] (ts)^{i}(1+t^2s^2) + \sum_{i\geq 0} [V_{i \omega_1+\omega_{j+1}}] (ts)^{i+1} +\\
+ \sum_{i\geq 0} [V_{i \omega_1+\omega_{j-1}}] (ts)^{i+1}.
\end{multline*}
We want to compute $P_{\SV{}{}}(t,s)$:
\begin{equation} \label{eq:P_intermedia}
P_{\SV{}{}}(t,s) = \sum_{\substack{1 \leq j \leq g \\i\geq 0 }}[\SV{j}{i}] t^{j+i}s^i = \sum_{\substack{0 \leq j \leq g \\i\geq 0 }} [V_{\omega_j}\otimes \SSS^i V] \frac{t^{2(g-j+1)}-1}{t^2-1} t^{j+i}s^i.
\end{equation}
In this last sum,  the addendum for $j=0$ is:
\[\sum_{i\geq 0} [V_{i\omega_1}] \frac{t^{2(g+1)}-1}{t^2-1} t^is^i = \frac{t^{2(g+1)}-1}{t^2-1} + \sum_{i \geq 0}[V_{i \omega_1 + \omega_1}] \frac{t^{2(g+1)}-1}{t^2-1} t^{i+1}s^{i+1}\]
and the addenda for $j>0$ are:
\begin{gather*}
\sum_{\substack{1 \leq j \leq g \\i\geq 0 }} [V_{i \omega_1+\omega_j}] \frac{t^{2(g-j+1)}-1}{t^2-1}  t^{j+i}s^{i}(1+t^2s^2) \\
 + \sum_{\substack{1 \leq j \leq g-1 \\i\geq 0 }} [V_{i \omega_1+\omega_{j+1}}] \frac{t^{2(g-j+1)}-1}{t^2-1}  t^{j+i+1} s^{i+1}\\
  + \sum_{\substack{1 \leq j \leq g \\i\geq 0 }} [V_{i \omega_1+\omega_{j-1}}] \frac{t^{2(g-j+1)}-1}{t^2-1}  t^{j+i+1} s^{i+1}.
\end{gather*}
Thus, eq.~\eqref{eq:P_intermedia} is equal to 
 \[ \begin{split}
\frac{t^{2(g+1)}-1}{t^2-1} + \frac{t^{2g}-1}{t^2-1} t^2s +
&\sum_{\substack{1 \leq j \leq g \\i\geq 0 }} [V_{i \omega_1+\omega_j}] \Bigl( \frac{t^{2(g-j+1)}-1}{t^2-1}  t^{j+i}s^{i}(1+t^2s^2) +\\
& + \frac{t^{2(g-j)}-1}{t^2-1} t^{j+i+2}s^{i+1} +\frac{t^{2(g-j+2)}-1}{t^2-1} t^{j+i}s^{i+1} \Bigr).
\end{split} \] 
Since the last factor is equal to $ (1+s)(1+t^2s) \frac{t^{2(g-j+1)}-1}{t^2-1}t^{i+j}s^i$, we have proven the claimed equality.
\end{proof}

\begin{lemma}\label{lemma:dim_rep}
For $i\geq 0$ and $1\leq j \leq g$, we have
\begin{equation}
\label{eq:dim_rep}
\dim \V{i}{j}= \binom{2g+i+1}{i, j 
} \frac{2g+2-2j}{2g+2+i-j} \frac{j}{i+j}.
\end{equation}
\end{lemma}
\begin{proof}
Recall that the positive roots of the Lie algebra $\spg$ are $e_k\pm e_h$ for $1\leq k<h\leq g$ and $2e_k$ for $1\leq k\leq g$.
Moreover, the half-sum of the positive roots is $\rho= \sum_{k=1}^g (g+1-k)e_k$.
Now, we apply the Weyl dimension formula:
\begin{align*}
&\prod_{k<h} \frac{\langle i\omega_1+\omega_j+\rho, e_k-e_h \rangle}{\langle \rho, e_k-e_h \rangle}=\frac{(g+i)!}{i!(g-j)! (j-1)!(i+j)}\\
&\prod_{k<h} \frac{\langle i\omega_1+\omega_j+\rho, e_k+e_h \rangle}{\langle \rho, e_k+e_h \rangle} = \frac{(2g+i+1)!(g+1-j)!(2g+2-2j)} {(g+i+1)! (2g+1-j)!(2g+i+2-j)}\\
&\prod_{k=1}^g \frac{\langle i\omega_1+\omega_j+\rho, 2e_k \rangle}{\langle \rho, 2e_k \rangle}=\frac{g+1+i}{g+1-j}.
\end{align*}
We obtain eq.~\eqref{eq:dim_rep} by multiplying the right hand sides of the above identities.
\end{proof}

\section{The cohomology of configuration spaces}

The case of the sphere ($g=0$) is essentially different from the case $g>0$ and our approach is useless since $\spg$ is trivial for $g=0$.
We refer to \cite{Sevryuk} for the following theorem.
\begin{theorem}
The rational homology of $\Uconf{n}(S^2)$ is:
\begin{align*}
& H^0(\Uconf{n}(S^2);\Q)=\Q && \textnormal{for all } n,\\
& H^2(\Uconf{1}(S^2);\Q)=\Q &&\\
& H^3(\Uconf{n}(S^2);\Q)=\Q && \textnormal{for } n\geq 3,\\
& H^k(\Uconf{n}(S^2);\Q)=0 && \textnormal{otherwise.}
\end{align*}
\end{theorem}

From now on we assume $g>0$ and so $\omega= \sum_{i=1}^g a_i b_i \neq 0$.
\begin{lemma}\label{lemma:gr_F and cohomology}
For $g>0$ the filtration $F_n A$ is strictly compatible with the differential.
Therefore, $\gr^{F_\bigcdot} H^\bigcdot (A,\dd) \simeq  H^\bigcdot (\gr^{F_\bigcdot} A, \gr^{F_\bigcdot} \dd)$.
\end{lemma}

\begin{proof}
We need to prove that $\Im \dd \cap F_n A \subseteq \dd( F_nA)$ for all $n\geq 0$.
Consider a generic element $x s1 +y p s1 + z + w p$ in $F_n A$ with $x,y,z,w \in \SV{}{}$.
Since the filtration $F_nA$ is induced by $\deg_3$, we can assume that $x \in F_{n-2}A$, $y \in F_{n-3}A$, $z \in F_n A$, and $w \in F_{n-1}A$.
Suppose that $\dd (x s1+y p s1+z+w p) \in F_{n-1} A$, then we have
 \[
 \dd(x s1+y p s1+z+w p)=  \tdd(x) p s1 + (-1)^{|x|} x (p -\omega ) - (-1)^{|y|} y p \omega + \tdd (z) p.
\] 
By looking at the third degree of the element 
in the right hand side,
 it follows that $\tdd(x) \in F_{n-4}A$, $\tdd(z)-\omega y \in F_{n-2}A$ and $\omega x \in F_{n-1}A$.
Since on $\SV{\bigcdot}{\bigcdot}$ the third grading coincides with the the total degree (i.e. $\deg_3(x)=|x|$), we can suppose $x,y,z$ being homogeneous of total degree $n-2$, $n-3$, and $n$ respectively.
So we have $\tdd(x)=0$, $\tdd(z)-\omega y=0$, and $\omega x=0$.
From $\omega x=0$ and $\tdd(x)=0$ we deduce that $\deg(x)>0$ and $x=\tdd(x')$ for some $x'$ of total degree $n-1$.
It follows that $\dd (x s1+y p s1+z+w p)=\dd (x')$ for $x' \in F_{n-1}A$ and so $\Im \dd \cap F_{n-1}A \subseteq \dd(F_{n-1}A)$.
\end{proof}

From now on we will work in $\gr_F A$ with the differential $\gr_F \dd$.
The only difference between $\dd$ and $\gr_F \dd$ is that $(\gr_F \dd)(s1)= -\omega$.
By an abuse of notation we denote the differential of $\gr_F A = A$ by $\dd$.

For any bigraded vector space $W$ we denote by $W[a,b]$ the same vector space with the bigraded shifted by $(a,b)$.

\begin{lemma}\label{lemma:ker}
The kernel of the differential $\dd$ is the direct sum of the following vector spaces:
\begin{enumerate}
\item $(\ker \tdd \cap \ker \Lo) [0,1]$,
\item $(\Im \tdd \cap \ker \Lo)[2,2] \oplus \ker \tdd [2,1]$, 
\item $\ker \tdd$,
\item $\SV{\bigcdot}{\bigcdot}[2,0]$.
\end{enumerate}
\end{lemma}

\begin{proof}
Consider a generic element $xs1+yps1+z+wp$ with $x,y,z,v \in \SV{\bigcdot}{\bigcdot}$: its differential is
\begin{equation}\label{eq:diff_generic}
 \dd(xs1+yps1+z+vp)=  \tdd(x) p s1 - (-1)^{|x|} x \omega - (-1)^{|y|} y p \omega + \tdd (z)p.
\end{equation}
Therefore $\dd(xs1+yps1+z+wp)=0$ if and only if $\omega x=0$, $\tilde{\dd}(x)=0$ and $(-1)^{|y|} y \omega = \tilde{\dd}(z)$.
The equations $\omega x=0$ and $\tilde{\dd}(x)=0$, imply that $x \in (\ker \tdd \cap \ker \Lo) $.
The condition $y \omega \in \Im \tdd$ is equivalent to $\tdd (y \omega)=\tdd(y) \omega =0$, thus 
\[y \in \ker \tdd \oplus (\Im \tdd \cap \ker \Lo)[0,1].\]
Let $z'=z'(y)$ be a fixed element such that $\omega y = \tdd (z')$:
then $z$ is of the form $z'+z''$ for some $z''\in \ker \tdd$ and $v$ can be any element in $\SV{\bigcdot}{\bigcdot}$.
\end{proof}

\begin{lemma}\label{lemma:im}
The image of the differential $\dd$ is the direct sum of the following vector spaces:
\begin{enumerate}
\item $0$,
\item $\Im \tdd [2,1]$,
\item $\omega \ker \tdd $,
\item $\Im \Lo [2,0] + \Im \tdd [2,0]$.
\end{enumerate}
\end{lemma}

\begin{proof}
Eq.~\eqref{eq:diff_generic} implies that the image of $\dd$ has trivial intersection with the submodule $s1\SV{\bigcdot}{\bigcdot}$.
Consider $x$ such that $\tdd(x)\neq 0$, then the element $(-1)^{|x|} x\omega  + \tilde{\dd}(x) p s1$ gives the addendum $\Im \tdd [2,1]$.
Now suppose $\tdd (x)=0$ and $x\omega \neq 0$, 
then $x\omega$ is in the image and generates a submodule isomorphic to $\omega \ker \tdd$.

Finally, $\Im \dd \cap p \SV{\bigcdot}{\bigcdot}$ coincides with $\Im \Lo [2,0] + \Im \tdd[2,0]$ (in general this is not a direct sum).
\end{proof}

Let $\oomega \in \SV{1}{1}$ be the unique $\spg$-invariant element such that $\tdd(\oomega)=\omega$.
The following lemma is an immediate consequence of \Cref{lemma:ker,lemma:im}.
\begin{lemma}\label{lemma:cohomology_module}
The cohomology $H^\bigcdot (A,\dd)$ is generated by:
\begin{enumerate}
\item[1.] $xs1$ for $x \in \ker \tilde{\dd} \cap \ker \Lo$,
\item[2.1.] $ps1+\oomega$, 
\item[2.2.] $yps1+x$ if $\tdd y  \in \ker \tdd \cap \ker \Lo$ and $(-1)^{|y|}y\omega = \tdd x$,
\item[3.] $y$ for $y \in \ker \tdd/\omega \ker \tdd$, 
\item[4.] $yp$ for $y \in \SV{\bigcdot}{\bigcdot} / (\Im \tdd + \Im \Lo)$.
\end{enumerate}
\end{lemma}

\begin{lemma}
The cohomologies of $\ker \Lo$, $\Im \Lo$, and $\coker \Lo$ with respect to the differential $\tilde{\dd}$ are given by:
\begin{align}
& H^{0,0}(\coker \Lo) =\langle 1 \rangle, \label{eq:1}\\
& H^{1,1}(\coker \Lo) =\langle \oomega \rangle, \label{eq:oomega} \\
& H^{2,0}(\Im \Lo)=\langle \omega \rangle, \label{eq:omega} \\
& H^{g,i}(\coker \Lo) =H^{g+1,i-1}(\Im \Lo)=H^{g,i-2}(\ker \Lo) \simeq \V{(i-2)}{g}, \label{eq:diff}\\
& H^{j,i}(\coker \Lo) =H^{j+1,i-1}(\Im \Lo)=H^{j,i-2}(\ker \Lo)=0. \label{eq:easy}
\end{align}
\end{lemma}
\begin{proof}
Consider the two short exact sequences
\begin{gather*}
0 \to \ker \Lo \to \SV{\bigcdot}{\bigcdot} \to \Im \Lo [2,0] \to 0\\
0 \to \Im \Lo \to \SV{\bigcdot}{\bigcdot} \to \coker \Lo \to 0.
\end{gather*}
By \Cref{lemma:exactness_tilde_dd} 
\begin{align*}
&H^{j,i}(\coker \Lo) \simeq H^{j+1,i-1}(\Im \Lo) &&\textnormal{ for } (j,i) \neq (0,0) \\
&H^{j,i}(\Im \Lo) \simeq H^{j-1,i-1}(\ker \Lo) &&\textnormal{ for } (j,i) \neq (2,0).
\end{align*}
Eq.~\eqref{eq:1}, \eqref{eq:oomega} and \eqref{eq:omega} follow immediately from the long exact sequence in cohomology.
Since $(\ker L_\omega)^{j,i}=0$ for $j<g$ and $(\coker L_\omega)^{j,i}=0$ for $j>g$, we deduce eq.~\eqref{eq:easy}.
The only representation that can appear in 
\[H^{g,i}(\coker \Lo) \simeq H^{g,i-2}(\ker \Lo)\]
is $\V{(i-2)}{g}$.
It is easy to see that the subspace $\V{i}{g} \subset V_{\omega_g} \otimes S^iV$ is contained in $\ker \Lo \cap \ker \dd$, but cannot lie in $\dd(\ker \Lo)$ since $(\ker \Lo)^{g-1,i+1}=0$.
This proves eq.~\eqref{eq:diff}.
\end{proof}

\begin{lemma}\label{lemma:ker_Lo_ker_d}
The Hilbert–Poincaré series of $\ker \tdd \cap \ker \Lo$ is
 \[ P_{\ker \tdd \cap \ker \Lo}(t,s)= t^{2g} + (1+t^2s) \sum_{\substack{1 \leq j\leq g \\ i \geq 0}} [\V{i}{j}] t^{2g-j+i}s^i.\] 
\end{lemma}
\begin{proof}
Notice that $\ker \tdd \cap \ker \Lo = \ker (\tdd_{|\ker \Lo})$ and 
\begin{align*}
P_{\ker \Lo}(t,s) =& \sum_{\substack{1 \leq j<g \\ i \geq 0}} [\V{i}{j}] t^{2g-j+i}s^i (1+s)(1+t^2s)+\\
 & +t^{2g}+t^{2g}s+\sum_{i \geq 0} [\V{i}{g}] t^{g+i}s^i (1+t^2s+t^2s^2), \\
P_{H(\ker \Lo)}(t,s) =&  \sum_{i \geq 0} [\V{i}{g}] t^{g+i}s^i.
\end{align*}
Using the formula $(s+1)P_{\ker (\tdd_{|\ker \Lo})}(t,s)= P_{\ker \Lo}(t,s)+s P_{H(\ker \Lo)}(t,s)$ we obtain the claimed equality.
\end{proof}

\begin{lemma}\label{lemma:ker_d/omegaker_d}
The Hilbert–Poincaré series of $\ker \tdd / \omega \ker \tdd$ is
 \[ P_{\ker \tdd / \omega \ker \tdd}(t,s)= 1 + (1+t^2s) \sum_{\substack{1 \leq j\leq g \\ i \geq 0}} [\V{i}{j}] t^{j+i} s^i.\] 
\end{lemma}
\begin{proof}
Consider the exact sequence
 \[ 0 \to \ker \tdd \cap \ker \Lo \to \ker \tdd \xrightarrow{\Lo} \ker \tdd \to \faktor{\ker \tdd}{\omega\ker \tdd} \to 0. \] 
We have
 \[P_{\ker \tdd/\omega\ker \tdd}(t,s)= (1-t^2) P_{\ker \tdd}(t,s) + t^2 P_{\ker \tdd \cap \ker \Lo}(t,s),\] 
and from \Cref{lemma:ker_Lo_ker_d} we obtain the claimed equality.
\end{proof}

\begin{lemma}\label{lemma:SV/im_Lo+im_d}
The Hilbert–Poincaré series of $\SV{\bigcdot}{\bigcdot}/\Im \Lo+ \Im \tdd$ is
 \[P_{\SV{}{}/\Im \Lo+ \Im \tdd}(t,s)=(1+t^2s)\big( 1+s \sum_{ \substack{ 1\leq j \leq g\\ i\geq 0}} [\V{i}{j}] t^{j+i} s^{i} \big).\] 
\end{lemma}
\begin{proof}
Let $K$ be the quotient $\SV{\bigcdot}{\bigcdot}/\Im \Lo+ \Im \tdd$.
Consider the exact sequence 
 \[ 0 \to \Im \tdd \cap \Im \Lo \to \Im \tdd \oplus \Im \Lo \to \SV{\bigcdot}{\bigcdot} \to K \to 0 \] 
and observe that $\Im \tdd \cap \Im \Lo= \ker \tdd \cap \Im \Lo = \ker (\tdd_{|\Im \Lo})$.
We compute the series $P_{\ker (\tdd_{|\Im \Lo})}$ using the formula
 \[ (1+s)P_{\ker (\tdd_{|\Im \Lo})}(t,s)= P_{\Im \Lo}(t,s) +  s P_{H(\Im \Lo)}(t,s)\] 
applied to the bigraded complex $(\Im \Lo, \tdd_{|\Im \Lo})$.
Notice that $P_{\SV{}{}}(t,s)= (1+s)P_{\Im \tdd}+1$ by \Cref{lemma:exactness_tilde_dd}, so we obtain
\begin{align*}
(1+s)P_{K}& =(1+s)\big(P_{\SV{}{}}-P_{\Im \tdd}-P_{\Im \Lo}+P_{\ker (\tdd_{|\Im \Lo})}\big)\\
&= sP_{\SV{}{}}+1-sP_{\Im \Lo}+sP_{H(\Im \Lo)}\\
&= sP_{\coker \Lo}+1+sP_{H(\Im \Lo)}.
\end{align*}
The equalities
\begin{align*}
P_{H(\Im \Lo)}(t,s) &= t^2 +t^2s \sum_{i\geq 0} [\V{i}{g}] t^{g+i} s^{i},\\
P_{\coker \Lo} &= 1+t^2s+ (1+s)(1+t^2s)\sum_{\substack{1\leq j <g\\ i\geq 0}} [\V{i}{j}] t^{j+i} s^{i} +\\
&+ (1+s+t^2s^2)\sum_{i \geq 0} [\V{i}{g}]t^{g+i}
\end{align*}
complete the proof.
\end{proof}

\begin{theorem}
The Hilbert-Poincaré series $P_{H(A)}(t,s)\in R_g[[t,s]]$ of $H(A_g, \dd)$ is
 \[
(1+t^2s)(1+t^2+t^{2g}s)+(1+t^2s)^2 \sum_{\substack{1\leq j \leq g\\ i\geq 0}} [\V{i}{j}] t^{j+i} s^{i} (1+t^{2(g-j)}s).
\] 
\end{theorem}

\begin{proof}
By \Cref{lemma:cohomology_module}:
 \[P_{H(A)}= (s+t^2 s^2)P_{\ker \tdd \cap \ker \Lo}+t^2s
  +P_{\ker \tdd/\omega \ker \tdd} +t^2 P_{\faktor{\SV{}{}}{\Im \Lo+ \Im \tdd}}.\] 
The computations of \Cref{lemma:ker_Lo_ker_d,lemma:ker_d/omegaker_d,lemma:SV/im_Lo+im_d} complete the proof.
\end{proof}

Let $Q_{g}(t,s,u)$ be the following series in the Grothendieck ring of $\spg$:
 \[Q_{g}(t,s,u) \defeq \sum_{i,j,n} [\gr_{i+2j}^W H^{i+j}(\Uconf{n}(\Sigma_g))] t^i s^j u^n .\] 

\begin{theorem}\label{theorem:main_theorem}
If $g>0$, the polynomial $Q_{g}(t,s,u) \in R_g[[t,s,u]]$ is equal to
\begin{align*}
Q_{g}& (t,s,u)= \frac{1}{1-u}\Bigl( (1+t^2su^3)(1 + t^2u) + (1+t^2su^2) t^{2g}su^{2(g+1)} + \\
& + (1+t^2su^2)(1+t^2su^3)\sum_{\substack{1\leq j \leq g\\ i\geq 0}} [\V{i}{j}] t^{j+i} s^{i} u^{j+2i} (1+t^{2(g-j)}su^{2(g-j+1)}) \Bigr).
\end{align*}
\end{theorem}
\begin{proof}
Use \Cref{lemma:gr_F and cohomology} and notice that $Q_K(t,s,u)=P_K(tu,su)$ for any sub-quotient $K$ of $\SV{\bigcdot}{\bigcdot}$, thus:
\begin{equation}
\begin{aligned}
\label{eq:Q_P_in_proof}
(1-u)Q_{g}=& su^2 P_{\ker \tdd \cap \ker \Lo}+t^2su^3+ t^2s^2u^4 P_{\ker \tdd \cap \ker \Lo} +\\ 
&+ P_{\ker \tdd/\omega \ker \tdd} +t^2 u P_{\SV{}{}/\Im \Lo+ \Im \tdd}.
\end{aligned}
\end{equation}
\Cref{lemma:ker_Lo_ker_d,lemma:ker_d/omegaker_d,lemma:SV/im_Lo+im_d} complete the proof.
\end{proof}

\Cref{theorem:main_theorem} and \Cref{lemma:dim_rep} give a formula for the mixed Hodge numbers and for the Betti numbers of $\Uconfn(\Sigma_g)$.
We use this formula to give a different proof of the result in \cite{DCK}*{Corollary 4.9} about the polynomial growth of the Betti numbers, and to extend it to the mixed Hodge numbers.

\begin{corollary}
For $n>k \gg 0$, the weights $h$ that appear in $H^k (\Uconfn(\Sigma_g))$ are in the range 
\[\min \left\{1,g-1\right\} \leq 3k-2h \leq \max \left\{ g+2, 2g \right\} \textnormal{ and } h \geq k.\]
Moreover,  in this range the dimension of $\gr^W_h H^k(\Uconfn(\Sigma_g))$ is polynomial in $k$ of degree $2g-1$.
\end{corollary}
\begin{proof}
For $i\gg 0$ the dimension of $V_{i\omega_1 + \omega_j}$ is polynomial in $i$ of degree $2g-1$ (see \Cref{lemma:dim_rep}).
\Cref{theorem:main_theorem} implies  $\gr_h^W H^k (\Uconfn(\Sigma_g))=0$ for $h<k$ or outside the range $\min \left\{1,g-1\right\} \leq 3k-2h \leq \max \left\{ g+2, 2g \right\}$.
For such $h$ and $k$, $\gr_h^W H^k (\Uconfn(\Sigma_g)$ is the direct sum of at most $8$ irreducible representations $V_{i\omega_1 + \omega_j}$ for some $i,j$ such that $i=k+\mathcal{O}(1)$.
It follows that $\dim(\gr_h^W H^k (\Uconfn(\Sigma_g))\sim i^{2g-1} \sim k^{2g-1}$ has polynomial growth in $k$.
Notice that, for any fixed $k$ the weights $h$ that appear in $H^k (\Uconfn(\Sigma_g))$ are at most $g+1$.
The claim about the Betti numbers follows since they are the sum of $g+1$ positive numbers (i.e.\ $\dim(\gr_h^W H^k (\Uconfn(\Sigma_g))$) that grow polynomially in $k$ of degree $2g-1$.
\end{proof}

Notice that the growth of Betti numbers of the unordered configuration space of the torus ($g=1$) is polynomial in $n$ of degree $2k-2$.
Indeed in \cite{PagAsymptotic} it is proven that $\beta_k(\Confn(\Sigma_1)) \sim \binom{n}{2k-2}$.

The same techniques can be applied to compute the invariants of configuration spaces of algebraic surfaces with zero irregularity.

\subsection*{Comparison with \texorpdfstring{\cite{DCK}}{[DCK]}}
Our work is dual to the previous one by Drummond-Cole and Knudsen, we briefly compare the two articles.
The Chevalley-Eilenberg complex $CE(\mathfrak{g}_{\Sigma_g})$ (\cite{DCK}*{Definition 2.1}) is dual to our differential algebra $B_g$.
Indeed, \Cref{lemma:iso_phi_n} and \cite{FT05}*{Theorems 1.14} imply that $H^\bigcdot(F_n B_g,\dd ) \simeq H^\bigcdot( \Uconfn (\Sigma_g))$.
Dually, the main result of \cite{Knudsen}*{Theorem 1.1} asserts that for $M=\Sigma_g$
\[\bigoplus_{n\geq 0} H_\bigcdot(\Uconfn (\Sigma_g))= H_\bigcdot (CE(\mathfrak{g}_{\Sigma_g})).\]

In \cite{DCK}*{Lemma 5.1} is proven that a complex $\mathcal{Z}\oplus v \mathcal{Z} $ is a deformation retract of $CE(\mathfrak{g}_{\Sigma_g})$; the dual of  $\mathcal{Z}\oplus v \mathcal{Z} $ is the algebra $A_g$ and the analogous statement is given by \Cref{lemma:I_is_acyclic}.
The submodule $\mathcal{K}_g$ of \cite{DCK} is dual to our module $P_{\SV{}{}/\Im \Lo+ \Im \tdd}$ (see \Cref{lemma:ker}).
The last part of the proofs, in which the Poincaré polynomial of $\mathcal{K}_g$ (resp.\ of $\SV{}{}/\Im \Lo+ \Im \tdd$) is computed, are essentially different:
ours uses the representation theory of the symplectic group while theirs uses homotopy and auxiliary spaces $\mathcal{V}(g,n)$.

Drummond-Cole and Knudsen identified a stable range for $H_i(\Uconfn (\Sigma_g))$ for $n>i$ and it is known that $H_i(\Uconfn (\Sigma_g))=0$ for $n<i-1$. 
Hence the unstable (non-trivial) range is for $n=i$ and $n=i-1$.
We observe the same phenomenon: in eq.\ \eqref{eq:main_in_introduction} the total degree in $t$ and $s$ differs from the degree in $u$ by at most one.
More precisely, the unstable polynomial $\mathcal{P}_1(t)$ of \cite{DCK}*{Theorem 4.2} is the Poincaré polynomial of the summand $t^2 u \SV{}{}/\Im \Lo+ \Im \tdd$ of eq.\ \eqref{eq:Q_P_in_proof} (corresponding to $n=i-1$).
The other unstable polynomial $\mathcal{P}_0(t)$ of \cite{DCK}*{Theorem 4.2} is the Poincaré polynomial of the summands
\[t^2su^3+ t^2s^2u^4 P_{\ker \tdd \cap \ker \Lo} +
+ P_{\ker \tdd/\omega \ker \tdd} + t^2 u P_{\SV{}{}/\Im \Lo+ \Im \tdd}.\]
of eq.\ \eqref{eq:Q_P_in_proof} (corresponding to $n=i$).

Notice that their work provides a uniform treatment for any genus $g\geq 0$, but ours excludes the case $g=0$ because we need $\omega \neq 0$.

Taking the dimension of both sides of our main eq.\ \eqref{eq:main_in_introduction}, we obtain a new formula for the Poincaré polynomial that is essentially different from the one  given by cases in \cite{DCK}*{Corollary 4.5, 4.6 and 4.7}.

We do not have a direct proof that our formula coincides with the one given in \cite{DCK}, but both provide the Betti numbers of configuration spaces on $\Sigma_g$, and a computational check shows that they agree in a very big range.
The code is available on request.
The following example shows in few cases that our formula for Betti numbers agrees with the exceptional value given in \cite{DCK}*{Corollary 4.5, 4.6, 4.7}.
\begin{example}
Consider $g=2$ \ie the case of genus two surface, the first terms (with respect the total degree in $t$ and $s$) of eq.\ \eqref{eq:main_in_introduction} are:
\begin{multline} 
\label{eq:example}
\frac{1}{1-u}\Bigl( 1+ [V_{\omega_1}]tu + t^2u + [V_{\omega_2}]t^2u^2  + t^2su^3+ [V_{2\omega_1}]t^2su^3 + [V_{\omega_2}]t^2su^4 +\\ + [V_{\omega_1}]t^3su^3 + [V_{\omega_1}]t^3su^4  + [V_{\omega_1+\omega_2}]t^3su^4 + [V_{\omega_1}]t^3su^5 + \dots \Bigr).
\end{multline}
The dimension of the $\mathfrak{sp}_{4}$-representations involved can be computed using eq.\ \eqref{eq:dim_rep} and they are:
\begin{align*}
&\dim V_{\omega_1}=4, \\
&\dim V_{\omega_2}=5, \\
&\dim V_{2 \omega_1}=10, \\
&\dim V_{\omega_1+\omega_2}= 16.
\end{align*} 
We deduce $\dim H^i(\Uconfn(\Sigma_2))$ from eq.\ \eqref{eq:example} by setting $s=t$ and by considering the dimension of the coefficient of $t^i u^n$: 
\begin{itemize}
\item $\dim H^1(\Uconfn (\Sigma_2))=4$ for $n\geq 1$,
\item $\dim H^2(\Uconf{1} (\Sigma_2))=1$ and  $\dim H^2(\Uconfn (\Sigma_2))=6$ for $n\geq 2$,
\item $\dim H^3(\Uconf{2} (\Sigma_2))=0$, $\dim H^3(\Uconf{3} (\Sigma_2))=11$, and  $\dim H^3(\Uconfn (\Sigma_2))=16$ for $n\geq 4$,
\item $\dim H^4(\Uconf{3} (\Sigma_2))=4$, $\dim H^4(\Uconf{4} (\Sigma_2))=24$, and  $\dim H^4(\Uconfn (\Sigma_2))=28$ for $n\geq 5$.
\end{itemize}
These numbers coincide with the one provided in \cite{DCK}*{Corollary 4.5, 4.6, 4.7}.
\end{example}

\bibliography{Config_surfaces}{}

\end{document}